\documentclass[12pt]{amsart}

\usepackage{amsmath, amscd, graphicx, latexsym, hyperref, times, rlepsf}

\oddsidemargin 0in \theoremstyle{plain}
\newtheorem*{Main}{Theorem}
\newtheorem{Thm}{Theorem}[section]
\newtheorem{Lem}[Thm]{Lemma}
\newtheorem{Cor}[Thm]{Corollary}

\newtheorem{Def}[Thm]{Definition}
\newtheorem{Rem}[Thm]{Remark}

\newtheorem{quest}[Thm]{Question}

\newcommand{\bfz}{{\mathbb{Z}}}
\newcommand{\Z}{{\mathbb{Z}}}

\newcommand{\OB}{\mathfrak{{ob}}}

\def\cp{\hbox{${\mathbb C} P^2$}}

\def\bn{\operatorname{bn}}
\def\sg{\operatorname{sg}}
\def\n{\operatorname{n}}
\numberwithin{equation}{section}

\begin{document}
\title[]{Invariants of contact structures from open books}

\author{John B. Etnyre}

\author{Burak Ozbagci}

\address{School of Mathematics \\ Georgia Institute
of Technology \\  Atlanta  \\ Georgia}

\email{etnyre@math.gatech.edu}

\address{School of Mathematics \\ Georgia Institute
of Technology \\  Atlanta  \\ Georgia and Department of Mathematics \\ Ko\c{c} University \\
Istanbul \\ Turkey}

\email{bozbagci@math.gatech.edu and bozbagci@ku.edu.tr}


\date{\today}

\thanks{J.B.E. was partially supported by the NSF CAREER Grant DMS-0239600
and NSF Focused Research Grant FRG-024466.
B.O. was partially supported by the Turkish Academy of Sciences and by
the NSF Focused Research Grant FRG-024466.}

\begin{abstract}
In this note we define three invariants of contact structures in
terms of open books supporting the contact structures. These
invariants are the support genus (which is the minimal genus of a
page of a supporting open book for the contact structure), the
binding number (which is the minimal number of binding components of
a supporting open book for the contact structure with minimal genus
pages) and the norm (which is minus the maximal Euler
characteristic of a page of a supporting open book).
\end{abstract}

\maketitle 

\section{Introduction}
Emmanuel Giroux has recently found an amazing correspondence between open book decompositions of 3--manifolds and
contact structures, \cite{gi}.  Specifically, he has shown that any contact
structure $\xi$ on $M$ is related to an open book through the Thurston-Winkelnkemper construction \cite{tw}.
This breakthrough has provided the basis for a much greater understanding of contact
structures and 3--manifold topology.

In this note we define three new invariants of contact structures
using Giroux's correspondence. Let $(M,\xi)$ be a contact manifold.
The first invariant is called the {\em support genus} and is simply
the minimal possible genus for a page of an open book that supports
$\xi.$ This invariant is denoted by $\sg(\xi).$ In \cite{et1} it was
shown that $\sg(\xi)=0$ for any overtwisted contact structure and
$\sg(\xi)=0$ forces any symplectic filling of $(M,\xi)$ to have
intersection form that embeds in a negative definite diagonalizable
quadratic form. In \cite{oss} it was shown that $\sg(\xi)=0$ implies
the Heegaard-Floer contact invariant of $\xi$ is reducible. Thus it
is clear that the support genus is tied to subtle geometric
properties of the contact structure. We do not add much to the
understanding of support genus here.

This paper will concentrate on a second invariant, the {\em binding
number}. The binding number of $\xi$ is the minimal number of
binding components for an open book supporting $\xi$ and that has
pages of genus $\sg(\xi).$ We denote the binding number by
$\bn(\xi).$ Note that it is important to restrict the open books
considered in the definition of binding number to those that have
minimal genus pages, because otherwise the binding number would
always be one (any contact structure is supported by an open book
with connected binding). We make several computations of the binding
number.
\begin{Main}
The binding number of an overtwisted contact structure $\xi$ on a
3-manifold is bounded in terms of the Euler class of $\xi.$ In
particular, the infinitely many overtwisted contact structures on a
manifold with Euler class zero (or torsional) take on finitely many
binding numbers. Moreover, if $M$ is a rational homology sphere then
there is a universal bound, depending on $M,$ on the binding number
of any overtwisted contact structure on $M.$
\end{Main}
Through examples we also show that the binding number is really an invariant of the contact structure in that it is not determined
merely  by the topology of the manifold $M.$ We also show that in certain cases the binding number can place restrictions on the
topology of the contact structure.
\begin{Main}
If $\sg(\xi)=1$ and $\bn(\xi)=1$ then the Euler class of $\xi$ is zero.
\end{Main}
We give an example showing that there are contact structures with
$\sg(\xi)=1$ but having non-zero Euler class. Thus the restriction
on $\bn(\xi)$ is necessary. This theorem follows from the following
theorem.
\begin{Main}
If $\xi$ is supported by an open book having two or less components in its binding and hyper-elliptic monodromy then
the Euler class of $\xi$ is zero.
\end{Main}

The last invariant of contact structures we define is the {\em norm}. The norm of a contact structure is minus the maximal
Euler characteristic of a page of an open book supporting $\xi.$ We denote the norm by $\n(\xi).$ We do not know if
the norm of a contact structure is determined by the support genus and binding number. We are only able to establish
\[\min\{ 2\sg(\xi) + \bn(\xi) -2, 2\sg(\xi) +1\} \leq \n(\xi)\leq 2\sg(\xi) + \bn (\xi) -2.\]
Thus, when $\bn(\xi)\leq 3,$ we know $\n(\xi)= 2\sg(\xi) + \bn (\xi)
-2$ but in general the norm might be smaller.

In Section~\ref{openbook} we recall the basic definitions involving open book decompositions and Giroux's correspondence 
between them and contact structures. The following section discusses the definitions of the invariants of plane fields on 3-manifolds.
In Section~\ref{invtsec} we give the definition of support genus and binding number. The following two sections discuss the
computation of the binding number for planar and, respectively, elliptic open books. Section~\ref{normsec} contains the definition of
the norm of a contact structure and discusses its relation to the other two invariants. We end with several fundamental questions
concerning these new invariants. 

\section{Open book decompositions}\label{openbook}
Suppose that for an oriented link $L$ in a closed and oriented
3--manifold $M$ the complement $M\setminus L$ fibers over the circle
as $\pi \colon M \setminus L \to S^1$ such that $\pi^{-1}(\theta) =
\Sigma_\theta $ is the interior of a compact surface bounding $L$,
for all $\theta \in S^1$. Then $(L, \pi)$ is called an \emph{open
book decomposition} (or just an \emph{open book}) of $M$. For each
$\theta \in S^1$, the surface $\Sigma_\theta$ is called a
\emph{page}, while $L$ the \emph{binding} of the open book. The
monodromy of the fibration $\pi$ is defined as the diffeomorphism of
a fixed page which is given by the first return map of a flow that
is transverse to the pages and meridional near the binding. The
isotopy class of this diffeomorphism is independent of the chosen
flow and we will refer to that as the \emph{monodromy} of the open
book decomposition. An open book $(L, \pi)$ on a 3--manifold $M$ is
said to be \emph{isomorphic} to an open book $(L^\prime,
\pi^\prime)$ on a 3--manifold $M^\prime$, if there is a
diffeomorphism $f: (M,L) \to (M^\prime, L^\prime)$ such that
$\pi^\prime \circ f = \pi$ on $M \setminus L$. In other words, an
isomorphism of open books takes binding to binding and pages to
pages.

An open book can also be described as follows. First consider the
mapping torus $$\Sigma_\phi= [0,1]\times \Sigma/(1,x)\sim (0,
\phi(x))$$ where $\Sigma$ is a compact oriented surface with $r$
boundary components and $\phi$ is an element of the mapping class
group $\Gamma_\Sigma$ of $\Sigma$. Since $\phi$ is the identity map
on $\partial \Sigma$, the boundary $\partial \Sigma_\phi$ of the
mapping torus $\Sigma_\phi$ can be canonically identified with $r$
copies of $T^2 = S^1 \times S^1$, where the first $S^1$ factor is
identified with $[0,1] / (0\sim 1)$ and the second one comes from a
component of $\partial \Sigma$. Now we glue in $r$ copies of
$D^2\times S^1$ to cap off $\Sigma_\phi$ so that $\partial D^2$ is
identified with $S^1 = [0,1] / (0\sim 1)$ and the $S^1$ factor in
$D^2 \times S^1$ is identified with a boundary component of
$\partial \Sigma$. Thus we get a closed $3$-manifold $M= \Sigma_\phi
\cup_{r} D^2 \times S^1 $ equipped with an open book decomposition
whose binding is the union of the core circles in the $D^2 \times S^1$'s
that we glue to $\Sigma_\phi$ to obtain $Y$. In conclusion, an
element $\phi \in \Gamma_\Sigma$ determines a $3$-manifold together
with an ``abstract" open book decomposition on it. Notice that by
conjugating the monodromy $\phi$ of an open book on a 3-manifold $M$
by an element in $\Gamma_\Sigma$ we get an isomorphic open book on a
3-manifold $M^\prime$ which is diffeomorphic to $M$.

Suppose that an open book decomposition with page $\Sigma$ is specified by $\phi \in \Gamma_\Sigma$.
Attach a $1$-handle to the surface $\Sigma$ connecting two points on $\partial \Sigma$ to obtain a
new surface $\Sigma^{\prime}$. Let $\gamma$ be a closed curve in $\Sigma^{\prime}$ going over the
new $1$-handle exactly once. Define a new open book decomposition with $ \phi^\prime= \phi \circ
t_\gamma \in \Gamma_{\Sigma^{\prime}} $, where $t_\gamma$ denotes the right-handed Dehn twist along
$\gamma$. The resulting open book decomposition is called a \emph{positive stabilization} of the one
defined by $\phi$. If we use a left-handed Dehn twist instead then we call the result a
\emph{negative stabilization}. The inverse of the above process is
called \emph{positive} (\emph{negative}) \emph{destabilization}. Notice that although the resulting
monodromy depends on the chosen curve $\gamma$, the 3--manifold specified by $(\Sigma^\prime,
\phi^\prime)$ is diffeomorphic to the 3--manifold specified by $(\Sigma, \phi)$.

\subsection{ Basic topology of a 3--manifold given as an open book.} \label{topology}
Given an abstract description of an open book in a 3--manifold $M$
with page a compact oriented genus $g$ surface $\Sigma$ with $r$
boundary components and monodromy $\phi \in \Gamma_\Sigma$. We can
determine the basic topology of $M$ by calculating its fundamental
group and its first homology group as follows. Fix a point $p_j$ on
the $j$-th boundary component of $\partial \Sigma$, for all $j=1,
\ldots, r$. We will calculate the fundamental group of $M$ based at
$p_1$. Let $a_1,\ldots, a_g, b_1,\ldots,b_g, c_1, \ldots, c_r$ be
the standard generators of $\pi_1(\Sigma)$ based at $p_1$, where
$c_i$'s correspond to loops around the boundary components and let
$\theta_j$ denote a loop based at $p_j$ which is transverse to all
the pages of the open book, for $j= 1, \ldots, r$. Then a
presentation of the fundamental group of the mapping torus
$\Sigma_\phi$ can be given as
\[
\pi_1(\Sigma_\phi)= \langle a_i,b_i,c_j, \theta_1 \;  \vert \; \prod_{i=1}^g [a_i b_i]
\prod_{j=1}^r c_j, \theta_1 a_i \theta_1^{-1}\phi_*(a_i^{-1}), \theta_1 b_i \theta_1^{-1}
\phi_*(b_i^{-1}),  \theta_1 c_j \theta_1^{-1}\phi_*(c_j^{-1})   \rangle
\]
where $1\leq i \leq g, 1\leq j\leq r$, $[a_i, b_i]= a_ib_i
a_i^{-1}b_i^{-1}$, and $\phi_*$ denotes the action of $\phi$ on
$\pi_1(\Sigma)$. Now connect the base point $p_1$ to $p_j$ by an arc
$\sigma_j \subset \Sigma$ and observe that the loop $
\theta_1\sigma_j\theta_j^{-1} \phi_* (\sigma_j^{-1})$ bounds a disk
in $\Sigma_\phi$. When we cap off the boundary component of
$\Sigma_\phi$ carrying the base point $p_1$ by a $D^2\times S^1$,
the loop $\theta_1$ will clearly bound a disk in $D^2 \times S^1$
and therefore it will become trivial in the resulting fundamental
group. Similarly when we cap off the boundary component of
$\Sigma_\phi$ carrying the point $p_j$, the above relation becomes
$\sigma_j\phi_*(\sigma_j^{-1})=1$. As a consequence, by capping off
all the boundary components of $\Sigma_\phi$ by $D^2 \times S^1$'s,
we get the following presentation of the fundamental group of $M$:
\[
\pi_1(M) = \langle\; a_i, b_i, c_j \; \vert \; \prod_{i=1}^g [a_i
b_i] \prod_{j=1}^r c_j, \;a_i\phi_*(a_i^{-1}),\; b_i
\phi_*(b_i^{-1}) , \; \sigma_j\phi_*(\sigma_j^{-1}) \;\rangle,
\]
where $i=1,\ldots, n$ and $j=2,\ldots, r.$ By abelianizing
$\pi_1(M)$ we get a presentation of the first homology group of $M$
as:
$$H_1(M)= \langle\; a_i, b_i, c_j \; \vert \; \;a_i-\phi_*(a_i),\; b_i-\phi_*(b_i) , \; \sigma_j
-\phi_*(\sigma_j) \; \rangle,$$ where $\phi_*$ now denotes the action of $\phi$ on $H_1(\Sigma)$.
Note that $c_1=-(c_2+\cdots +c_r)$ in $H_1 (\Sigma)$.

\subsection{Open books and contact structures}
We will assume throughout this paper that a contact structure
$\xi=\ker \alpha$ is coorientable (i.e., $\alpha$ is a global
1--form) and positive (i.e., $\alpha \wedge d\alpha >0 $ ). A
contact structure $\xi$ on $M$ is said to be supported by an open
book decomposition $(L,\pi)$ of $M$ if $\xi =\ker \alpha$ for some
contact form $\alpha \in \Omega^1(M)$ such that $\alpha ( L)
> 0$ and $d \alpha > 0$ is on every page. Thurston and Winkelnkemper
\cite{tw} have shown that every open book supports a contact
structure. Recently, Giroux \cite{gi} has proven a converse.
Specifically he has shown that every contact 3--manifold admits a
supporting open book and two open books supporting  the same contact
structure admit a common positive stabilization. Moreover two
contact structures supported by the same open book are isotopic. We
refer the reader to \cite{et2} and \cite{ozst} for more on the
correspondence between open books and contact structures.

\section{Invariants of plane fields} \label{euler}

An oriented 2--plane field $\xi$ in a 3--manifold has an Euler class
$e(\xi)$ (i.e., the first Chern class $c_1( \xi)$) and a
3-dimensional invariant $d_3(\xi).$ When the second cohomology of
the manifold has 2-torsion, the Euler class can be refined to an
invariant $\Gamma(\xi).$ In this section we will describe a method
to calculate these invariants for the underlying 2--plane field of a
contact structure $\xi$ supported by a given open book, starting
from an explicit factorization of the monodromy of the open book
into Dehn twists. We begin by computing the rotation number of
Legendrian knots sitting on the page of a supporting open book.

\subsection{Rotation number of the Legendrian realization of a curve on a page.} \label{rotation}

\begin{figure}[ht]
  \begin{center}
     \includegraphics{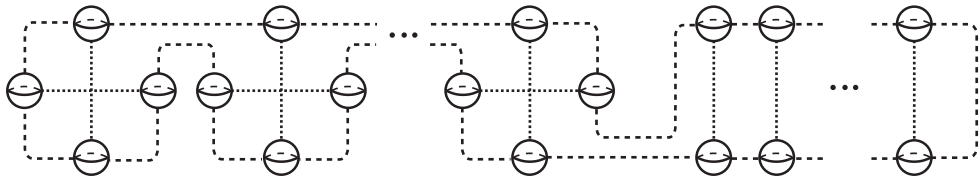}
   \caption{4--ball $D^4$ union $k$ 1--handles. The dashed line indicates the boundary of $\Sigma$ and
    dotted lines indicate generators of the homology of $\Sigma.$}
  \label{page}
    \end{center}
  \end{figure}

In this section we will compute the rotation number $r(\gamma)$ of
the Legendrian realization of a homologically nontrivial curve
$\gamma$ embedded on a page of an open book. We start with $P=
\natural_k S^1 \times D^3$, which is the 4--ball $D^4$ union $k$
1--handles, as shown in Figure~\ref{page}, where $k=2g+r$.
Evidently, $P$ is diffeomorphic to $D^2 \times \Sigma$, where
$\Sigma$ is a compact oriented surface of genus $g$ with $r$
boundary components. Note that $\partial P = \partial (D^2
\times\Sigma)= \partial D^2 \times \Sigma \cup D^2 \times \partial
\Sigma$ is a canonical decomposition of $\partial P$ into an open
book. The monodromy of this open book on $\partial P$ is the
identity map.

The way the attaching balls of the
1--handles are identified gives a trivialization of the tangent bundle of $P$. This trivialization also induces a
trivialization of the tangent bundle of the page which is depicted in Figure~\ref{page} with its
orientation induced from the blackboard: just take the usual oriented frame of $\mathbb{R}^2$,
restrict it to the disk and extend over the 1--handles.

Next we apply the construction of Thurston and Winkelnkemper in
\cite{tw} to put a contact structure $\xi$ on $\partial P =\#_k S^1
\times S^2$ supported by its canonical trivial open book
decomposition: Away from the binding the contact 1--form is given by
$\alpha =\beta + d\theta$, where $\theta$ is the meridional
direction to the binding, $\beta$ is a 1--form transverse to the
binding such that $d\beta$ is an area form on $\Sigma$. The pages of
the open book are convex and $\partial_{\theta}$ is the Reeb vector
field which is transverse to the pages.

Since the contact structure $\xi$ is supported by the trivial open
book on $\partial P$ by construction we can assume that the tangents
to the pages are as close as we wish to the contact planes, away
from the binding (cf. \cite{et2}). This gives a trivialization of
the contact planes restricted to a page. Recall that the rotation
number of a Legendrian curve is the winding number of its tangent
vector with respect to a fixed trivialization of the contact
structure over the curve. (Note the contact planes can be globally
trivialized since $e(\xi)=0$.) On the other hand a smooth curve on a
surface has a winding number which is defined as the winding number
of its tangent vector with respect to a fixed trivialization of the
tangent bundle of the surface. We conclude that the rotation number
of a Legendrian knot $\gamma$ is the same as the winding number of
the projection of $\gamma$ to a page $\Sigma$ since we can use the
same trivialization to define both invariants. This winding number
can be calculated for any curve $\gamma$ on the page, since we fix a
trivialization of the tangent bundle of the page once we fix the
attaching regions of the 1--handles in Figure~\ref{page}.

\subsection{Calculation of Euler class from
monodromy}\label{eclassmonod} Given an open book with page a genus g
surface $\Sigma$ with $r$ boundary components and monodromy $\phi$.
It is well-known that we can factor $\phi$ as a product of Dehn
twists along homologically nontrivial curves $\gamma_1, \gamma_2,
\cdots, \gamma_n$ on $\Sigma$. Below we will describe a method to
calculate the Euler class $e(\xi)$ (i.e., $c_1( \xi)$) of the
underlying 2--plane field of the contact structure $\xi$ supported
by this open book.

Note that an open book is the boundary of an achiral Lefschetz
fibration whose handlebody description as a 4--manifold is discussed
in detail in \cite{gs}. The curves $\gamma_1, \gamma_2, \cdots,
\gamma_n$ which appear in the factorization of the monodromy $\phi$
of the given open book are the vanishing cycles of this achiral
Lefschetz fibration. Consider $\gamma_i$'s as embedded in distinct
pages of the trivial open book on $\partial P = \#_k S^1 \times
S^2$, where $P= \natural_k S^1 \times D^3$, as depicted in
Figure~\ref{page}. We Legendrian realize each of these curves on
distinct pages. Then applying contact $(\pm 1)$--surgery on
$\gamma_i$ corresponds to adding a vanishing cycle to the achiral
Lefschetz fibration bounding the open book. If we apply
$(-1)$--surgeries only we get an ``honest'' Lefschetz fibration $X$
and the almost complex structure on $P$ will extend to $X$. The
Euler class of the contact structure $\xi$ is the restriction of
$c_1(X, J)$ to $\partial X=M$ and its Poincar\'{e} dual can be given
by :
\[
PD(e(\xi)) = \Sigma_{i=1}^n r(\gamma_i) [\mu_i] \in H_1(M,
\mathbb{Z})
\]
where $\mu_i$ is the meridian of $\gamma_i.$

If there exists $(+1)$--surgeries then although we can not find an
almost complex structure on the achiral Lefschetz fibration $X$,
there is an almost complex structure on $X \# q\cp$ (cf.
\cite{ozst}), where $q$ is the number of $(+1)$--surgeries. The
Poincar\'{e} dual $PD(e(\xi))$ will be evaluated in the same way as
above by restricting the first Chern class of the almost complex
structure on $X \# q\cp$ to the boundary.

When the second cohomology of $M$ has two torsion $e(\xi)$ is not the complete two dimensional
invariant of $\xi.$ In \cite{Gompf} a refinement of $e(\xi)$ was given. This invariant is a map
$\Gamma(\xi)$ from the spin structures on $M$ to $G=\{c\in H^2(M;Z)| 2c=e(\xi)\}.$ If $X$ is a Stein
2-handlebody obtained by attaching 2-handles to $D^4$ along a Legendrian link $L=\{K_1, \ldots,
K_n\}$ in $S^3=\partial D^4$ then $\Gamma$ can be described as follows: A spin structure
$\mathfrak{s}$ on $M=\partial X$ is described by a characteristic sub-link of $L_\mathfrak{s}\subset
L,$ see \cite{gs}. Then $\Gamma(\xi)(\mathfrak{s})$ is the restriction to $M$ of the
class $\rho\in H^2(X;Z)$ determined by
\begin{equation}\label{d2}
\langle \rho, \alpha_i\rangle = \frac 12 (r(K_i)+ \text{lk}(K_i, L_{\mathfrak{s}})),
\end{equation}
where $\alpha_i$ is the homology class in $X$ determined by $K_i.$

\subsection{The 3-dimensional invariant}
The 3-dimensional invariant $d_3(\xi)$ of a plane field is a rational number well-defined modulo the
divisibility of $e(\xi).$ We will only describe how to compute $d_3(\xi)$ when $e(\xi)$ is a torsion
element and thus $d_3(\xi)$ is well-defined. Let $M$ and $X$ be as above. Then we have
\begin{equation}\label{d3}
d_3(\xi)= \frac 14 (c^2(X) -3\sigma(X) - 2\chi(X)) + q,
\end{equation}
where $\sigma$ is the signature of $X,$ $\chi$ is the Euler
characteristic, and $q$ is the number of $(+1)$--surgeries. The
number $c^2(X)$ is the square of the class $c(X)$ with Poincar\'e
dual
\[
\sum_{i=1}^k r(\gamma_i) C_i,
\]
where the $C_i$'s are the cocores of the 2-handles attached along $\gamma_i$'s. Note that
$c(X)|_M=e(\xi)$, which we are assuming to be a torsion class. Thus some multiple $k\, c(X)$ of
$c(X)$, which naturally lives in $H^2(X; \bfz)$, comes from a class $c_r(X)$ in $H^2(X,\partial X; \bfz)$
which can be squared. So $c^2(X)$ means $\frac 1{k^2} c_r^2(X).$ Formula~\eqref{d3} is a
slight generalization of the one given in \cite{DingGeigesStipsicz04}, where it was assumed that $X$
had no 1-handles. Their proof caries over to our case, see \cite{EtnyreFuller}.

\section{Invariants of contact structures from open books}\label{invtsec}
There are several obvious invariants one can define using open book decompositions associated
to contact structures. We begin with the support genus.
\begin{Def}
The support genus of a contact structure $\xi$ on a 3-manifold $M$ is the minimal genus of a page of
an open book decomposition of $M$ supporting $\xi,$
\[
\sg(\xi)=\min \{\text{genus}(\Sigma) \vert (\Sigma, \phi) \text{ an open book decomposition supporting } \xi\}.
\]
\end{Def}
This definition was implicitly given and studied in \cite{et1} where the following was shown.
\begin{Thm}
If $(M,\xi)$ is a fillable contact structure and $\sg(\xi)=0$ then any filling of $(M,\xi)$
(1) has only one boundary component, (2) has negative definite intersection form and (3)
the intersection form embeds in a diagonalizable form.
\end{Thm}
\begin{Thm}\label{sgofot}
If $(M,\xi)$ is overtwisted then $\sg(\xi)=0.$
\end{Thm}
There are many examples of contact structures on lens spaces and
Seifert fibered spaces that have genus zero \cite{et1, Stephan}.
There are also many examples of contact structures having genus one
\cite{et1, EO}. Surprisingly, though there are many potential
examples \cite{EO}, it is still unknown if $sg(\xi)>1$ for any
contact structure $\xi.$

Next we define the binding number of a contact structure.
\begin{Def}
The binding number of a contact structure $\xi$ on a 3-manifold $M$ is the minimal number of components in the binding of the open book decomposition supporting $\xi$ that  have minimal genus,
\[
\begin{aligned}
\bn(\xi)=\{|\partial \Sigma| \vert (\Sigma,\phi) &\text{ an open book decomposition supporting }\\
& \xi \text{ with genus} (\Sigma)=sg(\xi)\}
\end{aligned}
\]
\end{Def}

\begin{figure}[ht]
  \relabelbox \small {
  \centerline{\epsfbox{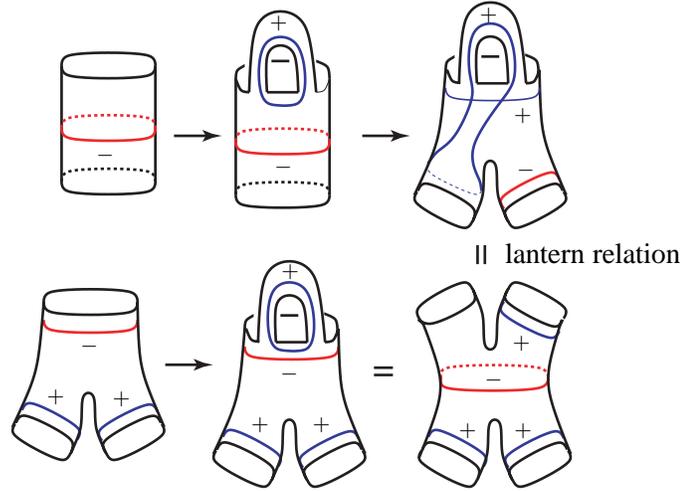}}}
  \relabel{1}{\tiny{$-$}}
  \relabel{2}{\tiny{$-$}}
  \relabel{3}{\tiny{$-$}}
  \relabel{4}{\tiny{$-$}}
  \relabel{5}{\tiny{$-$}}
  \relabel{0}{\tiny{$-$}}
  \relabel{6}{\tiny{$+$}}
  \relabel{7}{\tiny{$+$}}
  \relabel{8}{\tiny{$+$}}
  \relabel{9}{\tiny{$+$}}
  \relabel{a}{\tiny{$+$}}
  \relabel{b}{\tiny{$+$}}
  \relabel{c}{\tiny{$+$}}
  \relabel{d}{\tiny{$+$}}
  \relabel{e}{\tiny{$+$}}
  \relabel{f}{\tiny{$+$}}
  \relabel{g}{\tiny{$+$}}
  \relabel{l}{lantern relation}
  \endrelabelbox
        \caption{Positive stabilization of $\OB_1$ twice (the top row)
and $\OB_2$ once (the bottom row).}
        \label{iso}
\end{figure}

To illustrate the difficulty in computing the binding number of a
contact structure we give a simple example. Consider the contact
structures $\xi_1$ and $\xi_2$ in $S^3$ described by their contact
surgery diagrams depicted in Figure~\ref{iso}. The contact
structure $\xi_1$ can be given equivalently by two contact
$(+1)$--surgeries performed on the same unknot with its Legendrian
push-off. By computing the 3--dimensional invariants  of $\xi_i$ we
can see that $\xi_1$ is isotopic to $\xi_2$ and both are
overtwisted.

   \begin{figure}[ht]
  \relabelbox \small {
  \centerline{\epsfbox{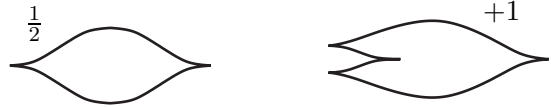}}}
  \relabel{1}{{$+1$}}
  \relabel{2}{{$\frac 12$}}
  \endrelabelbox
        \caption{Two isotopic overtwisted
contact structures $\xi_1$ (on the left) and $\xi_2$ (on the right)
in $S^3$.}
        \label{isotopic}
\end{figure}

The open books $\OB_1, \OB_2$ supporting  $\xi_1, \xi_2$,
are given in Figure~\ref{iso}, on the upper left and lower left, respectively. The $\pm$ sign on a
curve indicate a Dehn twist along that curve---right-handed for plus
and left-handed for minus sign. Thus it appears that $\bn(\xi_2)=3$
but in fact $\bn(\xi_2)=2.$ Using Giroux's correspondence we know that $\OB_1$ and $\OB_2$ are equivalent up to stabilization.
This is easily seen. If we positively
stabilize $\OB_1$ twice, and  $\OB_2$ once then
we can see that the resulting open books are
isomorphic  by the lantern relation on the four-holed sphere
as illustrated in Figure~\ref{iso}.

\section{Planar open books}
We first consider the binding number for contact structures
supported by planar open books. Our main goal in this section is to
give a bound on the binding number for overtwisted contact
structures and see that the binding number gives a non-trivial
invariant of a contact structure, that is that the binding number of
a contact structure can be larger than which is forced by the
topology of the manifold. We begin with a simple lemma.
\begin{Lem}\label{lem:basicplanar}
Suppose $\xi$ is a contact structure on a 3-manifold $M$ that is supported by a planar open book.
\begin{enumerate}
\item If $\bn (\xi)=1$ then $\xi$ is the standard tight contact structure on $S^3.$
\item If $\bn (\xi)=2$ and $\xi$ is tight then $\xi$ is the unique tight contact structure on $L(p,p-1)$ for some $p$.
\item If $\bn(\xi)=2$ and $\xi$ is overtwisted then $\xi$ is the overtwisted contact structure on $L(p,
1)$, for some $p$, with $e(\xi)=0$ and $d_3(\xi)=-\frac 14 p +\frac
34.$ When $p$ is even then the refinement of $e(\xi)$ is given by
$\Gamma(\xi)(\mathfrak{s})=\frac p2$ where $\mathfrak{s}$ is the
unique spin structure on $L(p,1)$ that extends over a two handle
attached to a $\mu$ with framing zero. Here we are thinking of
$L(p,1)$ as $-p$ surgery on an unknot and $\mu$ is the meridian to
the unknot.
\end{enumerate}
\end{Lem}
\begin{proof}
Statement (1) is obvious. An open book with annular pages supports a
tight contact structure if and only if the monodromy is $t_c^p$ for
some $p, $ where $t_c$ is a right handed Dehn twist about the core
circle $c$ in the annulus.  The contact structure supported by such
an open book can also be obtained by Legendrian surgery on $p-1$
copies of the Legendrian unknot with $tb=-1.$ This manifold is
easily seen be $L(p, p-1)$ and the contact structure must the unique
tight contact structure on it.

Similarly, the contact structure supported by an open book with
monodromy $t_c^{-p}$ can also be obtained by $(+1)$-contact surgery
on $p+1$ copies of the Legendrian unknot with $tb=-1.$ Thus the
manifold will be $L(p, 1)$ and the contact structure will be
overtwisted. The invariants $\Gamma$ and $d_3$ are easily computed
from Equations~\eqref{d2} and \eqref{d3}.
\end{proof}

\begin{Thm}
Let $\xi$ be the contact structure on the lens space $L(4,1)$ obtained from contact
surgery on the Legendrian unknot with $tb=-3$ and $r=0.$ Then $\sg(\xi)=0$ and $\bn(\xi)=4.$
However there is an open book for $L(4,1)$ that has annular pages.
\end{Thm}
This theorem indicates that the binding number of a contact
structure is sensitive to the plane field $\xi$ and is not
determined by the topology of the manifold. However, it turns out
that an open book for $L(4,1)$ that supports a plane field homotopic
to $\xi$ with fewer than four boundary components does not exist. So
$\bn(\xi)$ is still determined by the topology of the plane field.
It is interesting to note that there is an overtwisted contact
structure $\xi'$ on $L(4,1)$ with $\bn(\xi')=2$ with
$e(\xi)=e(\xi')$ and $d_3(\xi)=d_3(\xi').$ However,
$\Gamma(\xi)\not=\Gamma(\xi').$
\begin{proof}
One may easily construct an open book decomposition for $\xi.$ See Figure~\ref{fig:L4ex}.
\begin{figure}[ht]
  \relabelbox \small {
  \centerline{\epsfbox{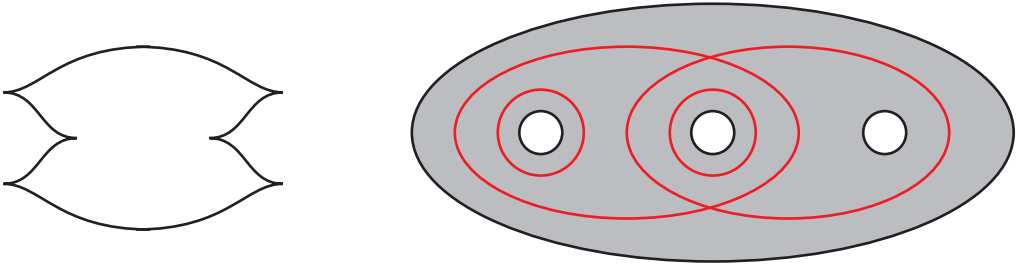}}}
  \endrelabelbox
        \caption{The contact structure $\xi$ on $L(4,1)$ is obtained by Legendrian surgery on the knot shown on the left.
            On the right is an open book supporting $\xi.$ The monodromy is a right handed Dehn twist about each curve shown.}
        \label{fig:L4ex}
\end{figure}
From this open book we see that $\bn(\xi)\leq 4.$ From
Lemma~\ref{lem:basicplanar} we see that $\bn(\xi)>2.$ Thus we must
rule out the possibility  $\bn(\xi)=3.$ To this end the surgery
picture for manifold given by open books with planar page having
three boundary components is given in Figure~\ref{fig:bn3}. For the
supported contact structure to be tight we must have $m, n$ and $k$
all be non-negative (otherwise, the open book will be non-right
veering \cite{hkm}). For this manifold to have a first homology
group of order 4, we must have $(m,n,k)=(0,4,1)$ or $(0,2,1).$ (Note
the surgery picture is symmetric in $m, n$ and $k$ so we do not list
all permutations.) One may easily check that $(m,n,k)=(0,4,1)$
yields $L(4,3)$ and $(m,n,k)=(0,2,2)$ yields $L(2,1)\# L(2,1)$ so
$\xi$ cannot be supported by a planar open book with three boundary
components.
\end{proof}

Recall that the set of homotopy classes of plane fields on a
3-manifold $M$ with a fixed Euler  class $e$ (or fixed refinement
$\Gamma$ if there is 2-torsion in $H^2(M;\Z)$) is in one-to-one
correspondence with $\Z_{d}$ where $d$ is the divisibility of $e$ in
$H^2(M;\Z).$ (Here $\Z_{d}$ will mean $\Z$ when $e$ is a torsion
element in $H^2(M;\Z)$.) Thus every manifold has infinitely many
homotopy classes of plane fields, and hence infinitely many
overtwisted contact structures. From \cite{et1} we know that the
support genus of any overtwisted contact structure is always zero.
It seems reasonable to believe that the binding number of an
overtwisted contact structure should be related to the homotopy
class of plane field. The next theorem indicates this relation is
very weak.
\begin{Thm}\label{thm:otb}
The binding number of an overtwisted contact structure $\xi$ on a
3-manifold is bounded in terms of the Euler class of $\xi.$ In
particular, the infinitely many overtwisted contact structures on a
manifold with Euler class zero (or torsional) take on finitely many
binding numbers.
\end{Thm}
\begin{Cor}
If $M$ is a rational homology sphere then there is a universal bound, depending on $M,$ on the binding number of any overtwisted
contact structure on $M.$
\end{Cor}
To prove this theorem we first study the overtwisted contact structures on $S^3.$ The range of $d_3$ for plane fields on $S^3$ is
$\Z + \frac 12,$ with the homotopy class containing the standard tight contact structure having $d_3=-\frac 12.$ We denote the unique
overtwisted contact structure on $S^3$ with three dimensional invariant $d_3$ by $\xi_{d_3}.$
\begin{Lem}\label{bnofotS3}
For any overtwisted contact structure $\xi$ on $S^3$ we have
\[
2\leq \bn(\xi)\leq 6.
\]
More precisely, we have
\[
\begin{aligned}
\bn(\xi_{\frac 12}) &= 2\\
\bn(\xi_{-\frac 12}) = \bn(\xi_{\frac 32}) &= 3\\
\bn(\xi_{\frac 52}) &= 4
\end{aligned}
\]
and
\[
\begin{aligned}
4\leq& \bn(\xi_{\frac {4p+1}{2}}) \leq 5, \qquad p\not = 0. \\
4\leq &\bn(\xi_{\frac{4p+3} 2}) \leq 6, \qquad p\not = -1,0.
\end{aligned}
\]
\end{Lem}
Note that this Lemma says Theorem~\ref{thm:otb} is true for $S^3.$
\begin{proof}
From Lemma~\ref{lem:basicplanar} we have $\bn(\xi_{\frac 12})=2$ and
$\bn >2$ for all other overtwisted contact structures on $S^3.$ We
next analyze the planar open books for $S^3$ with three binding
components. Let $\Sigma$ be the planar surface with three boundary
components. Any diffeomorphism of $\Sigma$ is determined by three
numbers $m,n,k,$ that give the number of Dehn twists on curves
$\gamma_1, \gamma_2, \gamma_3$ parallel to each boundary component.
Let $M_{m,n,k}$ be the 3-manifold determined by the open book with
page $\Sigma$ and monodromy given by
$t_{\gamma_1}^mt_{\gamma_2}^nt_{\gamma_3}^k.$ It is easy to see that
$M_{m,n,k}$ is the Seifert fibered space shown in
Figure~\ref{fig:bn3}.
\begin{figure}[ht]
  \relabelbox \small {
  \centerline{\epsfbox{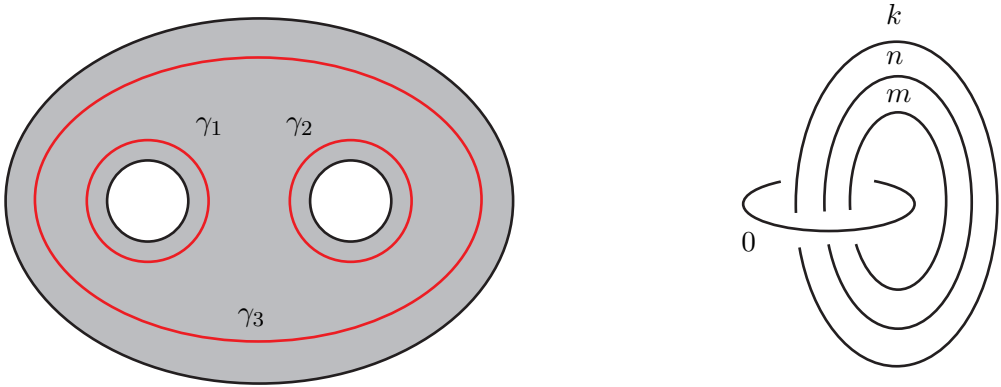}}}
  \relabel{1}{$\gamma_1$}
  \relabel{2}{$\gamma_2$}
  \relabel{3}{$\gamma_3$}
  \relabel{0}{$0$}
  \relabel{m}{$m$}
  \relabel {n}{$n$}
  \relabel {k}{$k$}
  \endrelabelbox
        \caption{The surface $\Sigma,$ left. The manifold $M_{m,n,k}$ right.}
        \label{fig:bn3}
\end{figure}
If $M_{m,n,k}$ is diffeomorphic to $S^3$ then we claim $|m|,|n|,
|k|$ cannot all be larger than one. Indeed, assume this is the case
and further assume that $m$ and $n$ are positive (this argument is
clearly symmetric in $m, n$ and $k$ so no generality is lost).
Figure~\ref{fig:sbn3}
\begin{figure}[ht]
  \relabelbox \small {
  \centerline{\epsfbox{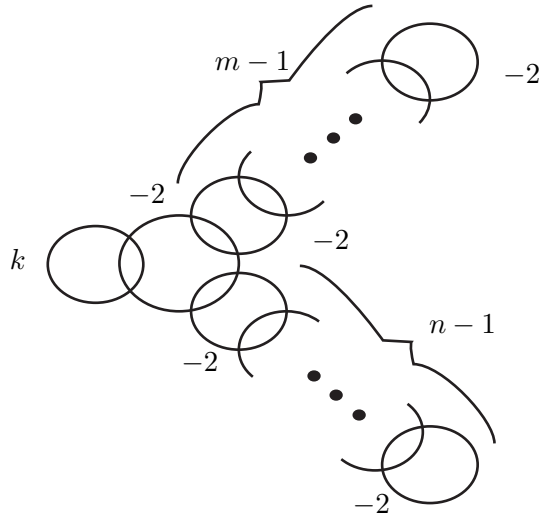}}}
  \relabel{1}{$-2$}
  \relabel{2}{$-2$}
  \relabel{3}{$-2$}
  \relabel{4}{$-2$}
  \relabel{5}{$-2$}
  \relabel {n}{$n-1$}
  \relabel {k}{$k$}
  \relabel{m}{$m-1$}
  \endrelabelbox
        \caption{Handle decomposition of a 4-manifold $X$ with $\partial X=M_{m,n,k}.$}
        \label{fig:sbn3}
\end{figure}
gives a handle decomposition of a 4-manifold $X$ with boundary
$M_{m,n,k}$. Moreover, if $k$ is negative, then this decomposition
can be realized by Legendrian knots with Thurston-Bennequin
invariants one larger than the framings. Thus $X$ has a Stein
structure. Since the only Stein filling $S^3$ has $b_2=0$ the
manifold $M_{m,n,k}$ cannot be $S^3.$ If $k$ is not negative then we
may expand the circle framed $k$ in to a leg of $k-1,$ $-2$--framed
unknots (as we did for the $m$ and $n$ framed unknots to go from
Figure~\ref{fig:bn3} to Figure~\ref{fig:sbn3}) changing the framing
on the central unknot to $-3.$ Finally, if $m$ and $n$ are less than
$-2$ we can reverse the orientation on $M_{m,n,k}$ and argue as
above.

To determine the contact structures realized form our open book for
$M_{m,n,k}$ we note that the corresponding contact surgery picture
is given in Figure~\ref{fig:csbn3}. (Recall $\frac 1m$ surgery on a
Legendrian knot is performed by taking $m$ Legendrian push-offs of
the Legendrian knot and performing $(+1)$, respectively $(-1)$,
contact surgery if $m$ is positive, respectively negative.)
\begin{figure}[ht]
  \relabelbox \small {
  \centerline{\epsfbox{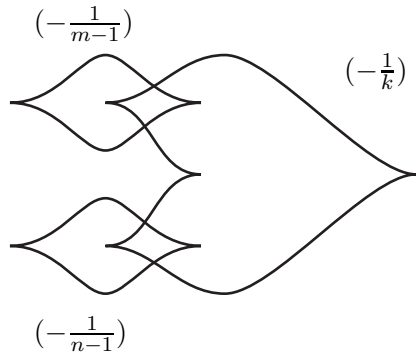}}}
  \relabel{1}{$(-\frac{1}{m-1})$}
  \relabel{2}{$(-\frac{1}{n-1})$}
  \relabel{3}{$(-\frac 1k)$}
  \endrelabelbox
        \caption{A contact surgery picture of $M_{m,n,k}$.}
        \label{fig:csbn3}
\end{figure}
We now know that $m,n$ or $k$ must be $0, \pm 1.$ If one of these is
zero then the other two must be $\pm 1.$ In this case we realize the
invariants $\frac 12$ (with $(m,n,k)=(-1,1,0)$), $\frac 32$ (with
$(m,n,k)=(-1,-1,0)$) and $-\frac 12$ (with $(m,n,k)=(0,1,1)$).
However, the contact structure with $d_3=-\frac 12$ is tight.

The other possibilities for $(m,n,k)$ are $(1,j,-1), (1,-1, j),
(-1,j,1),$ where $j$ is any integer, and $(3,2,-1), (3,-1,2),
(2,-1,3), (-2,-3,1), (-2,1,-3)$ and $(-3,1,-2).$ (Note interchanging
$m$ and $n$ does not affect the contact structure in
Figure~\ref{fig:csbn3} thus we do not list possibilities for
$(m,n,k)$ that differ by switching $m$ and $n.$ Actually, any
permutation of the number $(m,n,k)$ will lead to contactomorphic
contact structures, so one only really needs to check three cases.)
 For the first three sets of possibilities we get an overtwisted contact
structure with $d_3=\frac 12$ for then next three we get an overtwisted contact structure with $d_3=\frac 32$ and for the last three
we get overtwisted contact structures with $d_3=-\frac 12.$

Now for any overtwisted contact structure on $S^3$ with $d_3\not =
-\frac 12, \frac 12$ or $\frac 32$ we know that  $bn\geq 4.$ By
negatively stabilizing the open book for $\xi_{\frac 32}$ with three
boundary components we see that $d_3(\xi_{\frac 32})=4.$  For the
other overtwisted contact structures let $\xi_p$ be the contact
structure shown in  Figure~\ref{fig:aot}. It is easy to check that
$\xi_p$ is a contact structure on $S^3.$
\begin{figure}[ht]
  \relabelbox \small {
  \centerline{\epsfbox{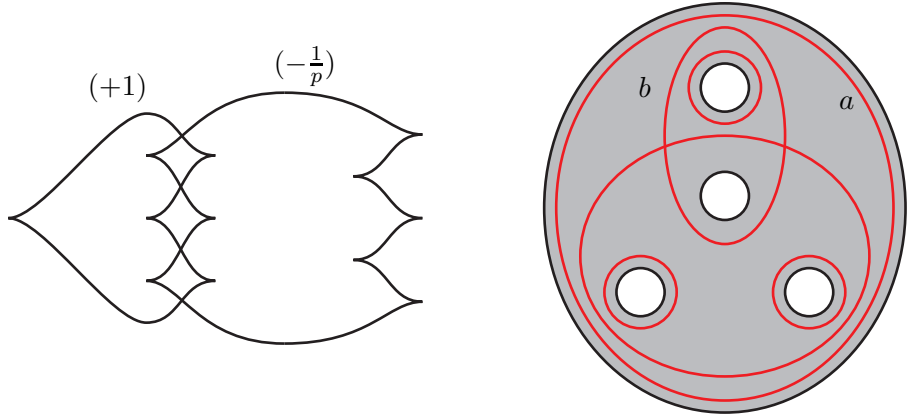}}}
  \relabel{1}{$(+1)$}
  \relabel{2}{$(-\frac{1}{p})$}
  \relabel{a}{$a$}
  \relabel{b}{$b$}
  \endrelabelbox
        \caption{A contact surgery picture of $\xi_p$ on the left. On the right is a surface $\Sigma.$}
        \label{fig:aot}
\end{figure}
Let $\Sigma$ be the surface shown on the right of
Figure~\ref{fig:aot} and $\phi=t_b^pt_a^{-1}\psi$, where $a$ and $b$
are the marked curves and $\psi$ is the diffeomorphism of $\Sigma$
consisting of a right-handed Dehn twist about the four unmarked
curves. The open book decomposition $(\Sigma,\phi)$ supports the
contact structure obtained by the contact surgery diagram on the
left in Figure~\ref{fig:aot}. From Figure~\ref{fig:aot} one may
easily compute
\[
d_3(\xi_p) = \frac{4p+1}{2}.
\]
Thus $\bn(\xi_{\frac{4p+1}{2}})\leq 5.$ If we negatively stabilize
$(\Sigma,\phi)$ we get an open book with planar pages for
$\xi_{\frac{4p+3}{2}}$ showing that $\bn( \xi_{\frac{4p+3}{2}})\leq
6.$
\end{proof}
\begin{Rem}{\em
When $(m,n,k)=(1,-1,j)$ in Figure~\ref{fig:csbn3}, one obtains the
same overtwisted contact structure on $S^3 $ independent of $j \in
\mathbb{Z}$. Thus we have an infinite family of open books $\OB_j$
for a single contact structure. According to Giroux any two in this
family will become isotopic after some number of positive
stabilizations. Using the lantern relation as we did in
Figure~\ref{iso} one can see that $j-j'$ positive stabilizations
will be sufficient to make $\OB_j$ and $\OB_{j'}$ isotopic.

On the other hand, the open books $\OB_p$ in Figure~\ref{fig:aot}
give an infinite family of distinct contact structures with
different $d_3$ invariants. Thus it is clear that $\OB_p$ must be
negatively stabilized $p'-p$ times to even give a contact structure
homotopic to the one supported by $\OB_{p'}$ where $p'>p.$ Thus for
the open books $\OB_p$ and $\OB_{p'}$ to give isotopic open books
$\OB_p$ must be stabilized at least $p'-p$ times. As a consequence
we get simple examples of fibered knots needing arbitrarily many
``Hopf plumbings" before they become isotopic. These examples are
simpler than those produced in \cite{nr}. }\end{Rem}

\begin{proof}[Proof of Theorem~\ref{thm:otb}]
Given two contact structures $\xi_1, \xi_2$ on $M^3$ that are
homotopic, as plane fields, over the two skeleton of $M$ there is an
overtwisted contact structure $\xi'$ on $S^3$ such that
$\xi_1\#\xi'$ is homotopic to $\xi_2$ as plane fields on all of $M.$
Thus given any overtwisted contact structure $\xi$ on $M$ we know
that it is supported by a planar open book with, say, $k$ boundary
components. By Murasugi summing the open book from
Lemma~\ref{bnofotS3} we get open book decompositions for all
overtwisted contact structures on $M$ with Euler class $e=e(\xi)$
having less than $k+7$ components in its binding. Hence we have a
bound for $\bn.$
\end{proof}

\section{Elliptic and hyper-elliptic open books}

We will call an open book whose page is genus one,
an elliptic open book.  An open book whose monodromy commutes with the hyper-elliptic involution will
be called hyper-elliptic.
In this section we consider restrictions on the Euler class of a contact structure that is supported by
an elliptic or hyper-elliptic open book.

\subsection{Elliptic open books}
From our discussion of rotations numbers in Section~\ref{rotation}
have the following simple observations. {\Lem \label{ellip} If $(M,
\xi)$ is supported by an elliptic open book with connected binding,
then $e(\xi)=0$. }

\begin{proof}
The monodromy $\phi$ of an elliptic open book with connected binding
can be expressed as a product of Dehn twists along the curves $a$
and $b$ depicted in Figure~\ref{elliptic}. But since the winding
numbers of $a$ and $b$ on the page are clearly both zero, the
rotation numbers of their Legendrian realizations are also zero, and
thus $e(\xi)=0$.
\begin{figure}[ht]
    \relabelbox \small {
  \centerline{\epsfbox{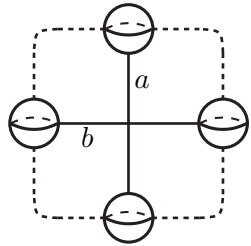}}}
  \relabel{a}{$a$}
   \relabel{b}{$b$}
  \endrelabelbox
\caption{Page of an elliptic open book with connected binding.}
     \label{elliptic}
\end{figure}

\end{proof}
It was conjectured that the Euler class of a fillable contact
structure supported by an elliptic open book is torsion. While the
lemma above verifies this conjecture if $\bn=1$ it is not always
true when the binding number is larger than one, as the following
example shows.
\begin{figure}[ht]
    \relabelbox \small {
  \centerline{\epsfbox{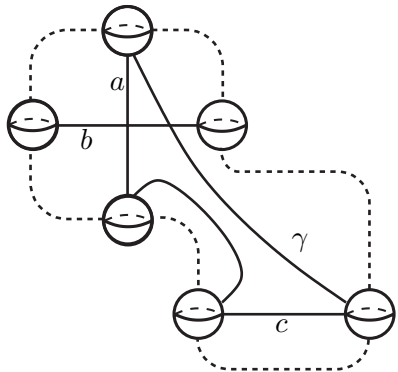}}}
  \relabel{a}{$a$}
   \relabel{b}{$b$}
   \relabel{c}{$c$}
   \relabel{g}{$\gamma$}
  \endrelabelbox
\caption{Page of an elliptic open book with two binding components.}
      \label{elltwo}
\end{figure}

Let $\Sigma_{1,2}$ denote a genus one surface with two boundary
components as in Figure~\ref{elltwo} and let $$\phi= (t_a t_b)^5
t_\gamma^2 t_c^2 : \Sigma_{1,2} \to \Sigma_{1,2}.$$ Let $(M,\xi)$ be
the contact 3-manifold given by the open book specified by
$(\Sigma_{1,2} , \phi).$
\begin{Thm}
The support genus and binding number of $(M,\xi)$ are
\[\sg(\xi)=1,\quad \bn(\xi)=2.\]
\end{Thm}

\begin{proof}
Be begin by computing the homology of $M$ using the presentation in
Section~\ref{topology}. The generators of $H_1 (M) $ are $a,b$ and
$c$ and the relations are
$$a+b = \phi_*(a)= a, \; \; -a +2b+ 2\gamma = \phi_* (b) =b, \; \; 2c+2\gamma=
\phi_*(\sigma)- \sigma=0,$$ where $\sigma$ is a curve connecting the
two boundary components. Hence  we get $$H_1 (M) = \langle a,c \;
\vert \;  a=-2c \rangle = \langle \gamma \rangle= \mathbb{Z},$$
where $\gamma= a+c$. On the other hand, when we apply the algorithm
in Section~\ref{eclassmonod}, the Euler class $e(\xi)$ of the
supported Stein fillable contact structure $\xi$ is given by
$2\gamma$ since the rotation numbers (i.e., the winding numbers) of
$a,b$ and $c$ are zero while the rotation number of $\gamma$ is
equal to one. Clearly $e(\xi) =2\gamma$ is non-torsion in $H_1(M)$.
Note that $(M, \xi)$ can not be supported by an elliptic open book
with connected binding by Lemma~\ref{ellip}.

We claim that $(M, \xi)$ is not supported by a planar open book
either. Consider the elliptic surface $E(1)$ with nine disjoint
sections. The monodromy of the elliptic fibration $E(1) \to S^2$ can
be given by $(t_b t_a)^{6}$, where $a$ and $b$ denote the standard
generators of the first homology group of a fiber. By removing the
union of a section and a cusp fiber from $E(1)$ we get an elliptic
fibration on the 4--manifold $W$ with once punctured torus fibers
whose monodromy is $(t_bt_a)^{5}$. One can check that $\partial W$
is diffeomorphic to $\Sigma(2,3,5)$ by Kirby calculus (see, for
example, \cite{gs}). Thus there is an induced open book on
$\Sigma(2,3,5)$ with monodromy $(t_bt_a)^{5}$. Since the monodromy
of this open book is a product of right-handed Dehn twists only, the
contact structure supported by this open book is Stein fillable (cf.
\cite{gi}) and hence isotopic to the unique tight contact structure
$\xi^\prime$ on $\Sigma(2,3,5)$. Note that, by \cite{et1},
$\xi^\prime$ can not be supported by  a planar open book. By
positively stabilizing the open book with monodromy $(t_bt_a)^{5}$
we get an open book (still supporting $\xi^\prime$) with two binding
components whose page is $\Sigma_{1,2}$ and whose monodromy is given
by $(t_a t_b)^5 t_c$. We now observe that $(M, \xi)$ is obtained
from $(\Sigma(2,3,5), \xi^\prime)$ by Legendrian surgeries. Hence
(the proof of) Theorem~1.2 of \cite{et1} implies  that
 $(Y, \xi)$ is not supported by a planar open book either. This fact also
follows from Corollary 1.4 in \cite{oss},
since $e(\xi)$ is non-torsion and $c^+(\xi) \neq 0$ (because $\xi$ is Stein fillable).
\end{proof}

\subsection{Hyper-elliptic open books}
\begin{figure}[ht]
  \relabelbox \small {
  \centerline{\epsfbox{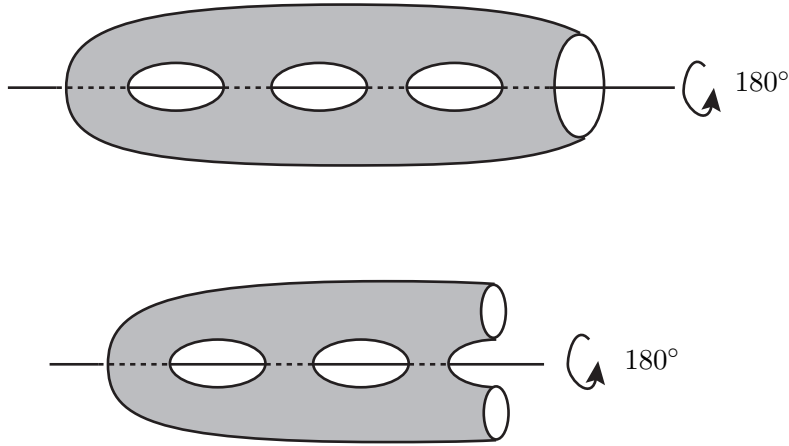}}}
  \relabel{1}{$180^\circ$}
  \relabel{2}{$180^\circ$}
  \endrelabelbox
        \caption{The hyper-elliptic involution $h_1$ top, and $h_2$ bottom.}
        \label{fig:he}
\end{figure}

Recall the hyper-elliptic involutions $h_1: \Sigma_{g,1}\to
\Sigma_{g,1}$ and $h_2:\Sigma_{g,2}\to \Sigma_{g,2}$ on a surface
with one or two boundary components, respectively, are the
involutions shown in Figure~\ref{fig:he}.
Note if the surface has genus $g$ then $h_1$ is the covering
translation for the unique 2--fold branch cover over $D^2$ with
$2g+1$ branch points and $h_2$ is the covering translation for the
unique 2--fold branch cover of $D^2$ with $2g$ branch points. We
call an open book $(\Sigma_{g,i}, \phi)$ {\em hyper-elliptic} if
$\phi\circ h_i= h_i\circ \phi,$ for $i=1,2.$ Lemma~\ref{ellip} above
is a special case of the following theorem.
\begin{Thm}\label{thm:he}
Let $(M,\xi)$ be supported by a hyper-elliptic open book. Then
$e(\xi)=0.$
\end{Thm}
\begin{Rem}{\em
We note the converse of this theorem is not true. Let $\xi$ be the
Stein fillable contact structure on $T^3.$ It is well known that
$e(\xi)=0$ but $\xi$ cannot be supported by a hyper-elliptic open
book since that would imply $T^3$ is a 2--fold branch cover over
$S^3$ (see the lemma below); but $T^3$ is not a 2--fold (or any
cyclic) branch cover over $S^3,$ \cite{hn}. }\end{Rem}

Before proving the theorem we recall the definition of a contact
branched cover, \cite{Geiges}.  Let $(M,\xi)$ be a contact
3-manifold and $\Gamma$ be a transverse link in it. Let $p:M'\to M$
be a covering map branched over $\Gamma.$ There is a standard
neighborhood $N$ of $\Gamma$ and $M'\setminus p^{-1}(N)$ is a
covering space of $M\setminus N. $ Thus we may pull $\xi$ back to
$M'\setminus p^{-1}(N).$ It is easy to extend this contact structure
over $M'$ so that $p^{-1}(\Gamma)$ is a transverse link and the
contact structure is the pull back of $\xi$ on $M'\setminus
p^{-1}(\Gamma).$ Call this contact structure $\xi'.$
\begin{Lem}
Suppose $(M,\xi)$ is supported by a hyper-elliptic open book $(\Sigma, \phi).$ Then $M$ is a 2--fold branch cover over $S^3$
and $\xi$ is the pull-back of the standard tight contact structure on $S^3.$
\end{Lem}
\begin{proof}
We consider the case when $\Sigma$ has one boundary component
leaving the other, analogous case to the reader. Recall $M$ is the
union of the mapping torus $T_\phi$ of $\phi$ and a solid torus $N.$
The mapping cylinder is obtained from $\Sigma\times [0,1]$ by
identifying $\Sigma\times\{0\}$ and $\Sigma\times\{1\}$ via $\phi.$
The hyper-elliptic involution $h_1$ of $\Sigma$ induces a branch
covering of $\Sigma\times [0,1]$ over $D^2\times [0,1]$ and since
$\phi$ commutes with $h_1$ we obtain a branch covering of $T_\phi$
over $D^2\times S^1.$ Note that this is an regular 2--fold coving of
$\partial T_\phi$ over $\partial D^2\times S^1.$ The solid torus $N$
double covers the solid torus $S^1\times D^2.$ Gluing $D^2\times
S^1$ and $S^1\times D^2$ together (in the obvious way that preserves
the product structure on the boundary) yields $S^3$ and the
coverings maps on $T_\phi$ and $N$ fit together to give a branch
covering of $M$ over $S^3$. It is also easy to see that the standard
open book with disk pages for $S^3$ pulls back under this covering
map to the open book for $M.$ Thus $\xi$ is the pull-back of the
standard contact structure on $S^3.$
\end{proof}

\begin{Lem}\label{lem:ec}
Let $(M,\xi)$ be a contact 3--manifold and $\Gamma$ a transverse link in it. Let $(M',\xi')$ be the
2-fold cover of $(M,\xi)$ over $\Gamma.$ Then
\[e(\xi')=p^*(e(\xi)- u),\]
where $u$ is the Poincar\'e dual of the homology class $[\Gamma].$
\end{Lem}
\begin{proof}
Let $s:M\to \xi$ be a section of $\xi$ transverse to the zero
section. Let $Z=s^{-1}(s).$ The Poincar\'e dual of $e(\xi)$ is
$[Z].$ Let $N$ be a neighborhood of $\Gamma$ as above. We can
homotop $s$ so that $Z$ is disjoint from $N.$ On $p^{-1}(M\setminus
N)$ $s$ will give a section $s'$ of $\xi'$ that is transverse to the
zero section and $(s')^{-1}(0)= p^{-1}(Z).$ We now extend this
section over $p^{-1}(N).$ To this end, note that each component of
$N$ is $S^1\times D^2$ and we may homotop $\xi$ in $N$ so that it is
tangent to $D^2$ factor. (We have left the world of contact
structures, but we are only trying to compute the Euler class of the
plane field $\xi'$ which can be done using any plane field homotopic
to it.) For each point $p\in \Gamma$ the vector $s(p)\in\xi$ is
transverse to $\Gamma.$ If $N$ is taken to be a sufficiently small
neighborhood of $\Gamma$ we may assume that $s(p)$ is the constant
vector field on each $\{pt\}\times D^2.$ Consider the branch map
$D^2\to D^2: z\to z^2.$  This map pulls back the constant vector
field along $\partial D^2$ to a non-zero vector field near $\partial
D^2$ with winding $-1.$ Thus it may be extended over $D^2$ to a
vector field with one transverse zero. Hence the section $s$ on
$\partial N$ induces a section of $\xi'$ on $\partial p^{-1}(N)$
that can be extended over $p^{-1}(N)$ so that it is zero along
$p^{-1}(\Gamma).$ We now have a section $s'$ of $\xi'$ defined on
all of $M'$ that is transverse to the zero section on $(s')^{-1}(0)=
p^{-1}(Z)-p^{-1}(\Gamma).$ Taking Poincar\'e duals we obtain the
desired formula.
\end{proof}

\begin{Rem}{\em
One may similarly derive a formula for more general branched covers. In particular, the
formula for other cyclic covers and simple coves is a direct generalization of the one
in Lemma~\ref{lem:ec}.
}\end{Rem}

\begin{proof}[Proof of Theorem~\ref{thm:he}]
If $(M,\xi)$ is a contact manifold given by the hyper-elliptic open
book $(\Sigma,\phi),$ then there is a 2--fold branched covering
$p:M\to S^3$ along $\Gamma\subset S^3$ such that
$\xi=p^*(\xi_{std}).$ There is a surface $S$ in $S^3$ such that
$\partial S = \Gamma.$ The 2--fold branched cover $M$ is formed by
gluing two copies of $S^3\setminus S$ together. Thus the branch set
$\Gamma$ in $M$ is the boundary of $S\subset M$ and hence
null-homologous. Lemma~\ref{lem:ec} now implies that $e(\xi)=p^*
e(\xi_{std})=p^* 0=0.$
\end{proof}

\section{The norm of a contact structure}\label{normsec}
In this section we discuss another invariant of contact structures derived form open book decompositions.
Motivated by the Thurston norm we define the norm of a contact structure.
\begin{Def}
The norm of a contact structure $\xi$ on a 3-manifold $M$  is minus the minimal Euler characteristic of a page of
an open book supporting $\xi,$
\[
\n(\xi)=\min \{ -\chi(\Sigma)| (\Sigma, \phi) \text{ an open book
decomposition supporting } \xi\}.
\]
\end{Def}
For a given surface any pair of ``genus", ``number of boundary
components" and ``Euler characteristic" determine the third.
However, it is not clear if the support genus and binding number of
a contact structure determine the norm. We have the following simple
observations.
\begin{Lem}\label{bound}
For any contact structure $\xi$ we have $\n(\xi)\leq -1$ with equality if and only if $\xi$ is the standard tight contact
structure on $S^3.$ Moreover, we have
\[\min\{ 2\sg(\xi) + \bn(\xi) -2, 2\sg(\xi) +1\} \leq \n(\xi)\leq 2\sg(\xi) + \bn (\xi) -2.\]
Thus for contact structures with $\bn(\xi)\leq 3$ we know $\n(\xi)=2\sg(\xi) + \bn (\xi) -2.$
\end{Lem}
\begin{proof}
The only non-trivial statement is the lower bound. If the minimum in the definition of $\n(\xi)$ is achieved with an
open book with genus $g>\sg(\xi)$ and $m$ boundary components then $\n(\xi)=2g+m-2.$ The smallest possible
value for the right hand side is $2\sg(\xi)+1.$
\end{proof}

\begin{Cor} For any contact structure $\xi$ in $S^3$ we have
$n(\xi) \leq 4$.
\end{Cor}

\begin{proof}
For any contact structure $\xi$ in $S^3$ we have $sg(\xi)=0$. The
corollary follows by combining Lemma~\ref{bnofotS3} with
Lemma~\ref{bound}.
\end{proof}

Similarly the next result follows from Theorem~\ref{thm:otb} and
Lemma~\ref{bound}.

\begin{Cor}
The norm of an overtwisted contact structure $\xi$ on a 3-manifold
is bounded in terms of the Euler class of $\xi.$ In particular, the
infinitely many overtwisted contact structures on a manifold with
Euler class zero (or torsional) take on finitely many norms.
\end{Cor}

\section{Questions}
From \cite{et1} and \cite{oss} we know the support genus has
geometric meaning in that if $\sg(\xi)=0$ this forces any filling of
$\xi$ to have certain properties and the Heegaard-Floer invariant to
have certain properties.
\begin{quest}
Is there geometric content to the binding number and the norm of a contact structure?
\end{quest}
It would be reasonable to suspect that the binding number has something to do with the homotopy class of plane field for the
contact structure, but Theorem~\ref{thm:otb} shows that is not the case or at least that the binding number is insensitive to
the homotopy class on the three skeleton. We ask
\begin{quest}
On a manifold $M$ is there a bound on the binding number of any overtwisted contact structure (that does not depend on
the Euler class, but only on $M$)?
\end{quest}
One could also ask if the binding number of tight contact structures is bounded. For a specific question consider
\begin{quest}
Let $\xi_n$ be the contact structure on $T^3$ with Giroux torsion $n.$ Is $\bn(\xi_n)$ bounded independent of $n$?
\end{quest}
Based on open book decompositions constructed by Van Horn \cite{vh}
one might conjecture: $\bn(\xi_n)=3n.$

Though there are countless other questions one could ask we end with
\begin{quest}
Is the norm of a contact structure determined by the support genus and binding number?
\end{quest}


\begin{thebibliography}{99999}



\bibitem[DGS]{DingGeigesStipsicz04}
F.~Ding, H.~Geiges, and A.~I.~Stipsicz,
{\em Surgery diagrams for contact 3-manifolds},
Turkish J. Math, {\bf 28} (2004) no.~1, 41--74.

\bibitem[E1]{et1}
J. B. Etnyre, {\em Planar open book decompositions and contact
structures,} IMRN {\bf 79} (2004), 4255--4267.

\bibitem[E2]{et2}
J. B. Etnyre, {\em Lectures on open book decompositions and contact
structures}, Lecture notes from the Clay Mathematics Institute
Summer School on Floer Homology, Gauge Theory, and Low Dimensional
Topology at the Alfr\'{e}d R\'{e}nyi Institute;
arXiv:math.SG/0409402.

\bibitem[EF]{EtnyreFuller}
J.~B.~Etnyre and T.~Fuller, {\em Realizing 4-manifolds as achiral Lefschetz fibrations}, to appear IMRN.

\bibitem[EO]{EO}
J.~B.~Etnyre and B.~Ozbagci, {\em Open books and plumbings}, preprint 2006.

\bibitem[Ge]{Geiges}
H.~Geiges,
{\em Constructions of contact manifolds},
Math. Proc. Cambridge Philos. Soc. {\bf 121} (1997), no. 3, 455--464.

\bibitem[Gi]{gi}
E. Giroux, \emph{G\'{e}ometrie de contact: de la dimension trois
vers les dimensions sup\'{e}rieures,} Proceedings of the
International Congress of Mathematicians (Beijing 2002), Vol. II,
405--414.

\bibitem[G]{Gompf}
R.~Gompf,
{\em Handlebody construction of {S}tein surfaces},
Ann. of Math. (2), {\bf 148} (1998) no.~2, 619--693.

\bibitem[GS]{gs}
R. Gompf and A. Stipsicz, {\em 4--manifolds and Kirby calculus},
Grad. Stud. Math., Vol. {\bf 20}, AMS, 1999.

\bibitem[HN]{hn}
U.~Hirsch and W.D.~Neumann, {\em On cyclic branched coverings of spheres},
Math. Ann. {\bf 215} (1975), 289--291.

\bibitem[HKM]{hkm}
K.~Honda, W.~Kazez and G.~Matic, {\em Right-veering diffeomorphisms of compact surfaces with boundary I}, preprint 2005.


\bibitem[NR]{nr}
W.D.~Neumann and L.~Rudolph, {\em Difference index of vector fields and the enhanced Milnor number,}
Topology {\bf 29} (1990) no.~1, 83--100.

\bibitem[OS]{ozst}
B. Ozbagci and A. Stipsicz, {\em Surgery on contact 3--manifolds and
Stein surfaces}, Bolyai Soc. Math. Stud., Vol. {\bf 13}, Springer,
2004.

\bibitem[OSS]{oss} P. Ozsvath, A. Stipsicz and Z. Szabo,
{\em Planar open books and Floer homology}, arXiv: math.SG/0504403

\bibitem[Ro]{rol}
D. Rolfsen, \emph{Knots and links,} Publish or perish, 1976.

\bibitem[S]{Stephan}
S.~Sch\"onenberger, {\em Planar open books and symplectic fillings}, PhD.~Dissertation, the University of
    Pennsylvania, 2005.

\bibitem[TW]{tw}
W. Thurston and H. Winkelnkemper, \emph{On the existence of contact forms,} Proc. Amer.
Math. Soc. \textbf{52} (1975), 345--347.

\bibitem[Vh]{vh}
J. Van Horn, PhD.~Dissertation, the University of Texas at Austin,
in preparation.

\end{thebibliography}
\end{document}